\documentclass[12pt]{amsart}

\usepackage{amsmath,amsfonts,amssymb,amsthm,mathtools}

\usepackage{verbatim}
\usepackage{amsmath,amsfonts}
\usepackage[mathscr]{euscript} 
\usepackage{amsthm}
\usepackage{blkarray,bigstrut} 
\usepackage{url}

\usepackage{graphicx} 

\usepackage{enumitem} 

\usepackage{hyperref} 
\hypersetup{colorlinks} 

\usepackage[colorinlistoftodos]{todonotes}

\RequirePackage{cleveref}
\usepackage{hypcap}
\hypersetup{colorlinks=true, citecolor=darkblue, linkcolor=darkblue}
\definecolor{darkblue}{rgb}{0.0,0,0.7}
\newcommand{\darkblue}{\color{darkblue}}

\definecolor{darkred}{rgb}{0.68,0,0}

\definecolor{darkgreen}{rgb}{0,.38,0}

\newcommand{\defnb}[1]{\emph{\darkblue #1}}

\setlist[enumerate]{
	label=\textnormal{({\roman*})},
	ref={\roman*}}

\makeatletter
\def\th@plain{%
	\thm@notefont{}
	\itshape 
}
\def\th@definition{%
	\thm@notefont{}
	\normalfont 
}
\makeatother

\newtheorem{thm}{Theorem}[section]
\newtheorem{lemma}[thm]{Lemma}

\newtheorem{prop}[thm]{Proposition}

\newtheorem{question}[thm]{Question}

\theoremstyle{definition}

\newtheorem{definition}[thm]{Definition}

\numberwithin{figure}{section}
\numberwithin{equation}{section}

\def\aa{\mathrm{a}}

\DeclareMathOperator{\ax}{\mathrm{x}} 
\DeclareMathOperator{\ay}{\mathrm{y}} 

\DeclareMathOperator{\ba}{\mathbf{a}} 
\DeclareMathOperator{\bx}{\mathbf{x}} 
\DeclareMathOperator{\by}{\mathbf{y}} 
\def\cA{\mathcal A}
\def\cB{\mathcal B}

\DeclareMathOperator{\indeg}{\mathrm{indeg}} %
\def\nn{\mathbb N}

\DeclareMathOperator{\occ}{\mathrm{occ}} 
\DeclareMathOperator{\outdeg}{\mathrm{outdeg}} %

\def\rr{\mathbb R}

\DeclareMathOperator{\bZ}{\mathbf{Z}} 
\DeclareMathOperator{\az}{\mathrm{z}} 
\DeclareMathOperator{\Rb}{\mathbb{R} } 
\DeclareMathOperator{\aM}{\mathrm{M}} 

\DeclareMathOperator{\Out}{\mathrm{Out}} 
\DeclareMathOperator{\In}{\mathrm{In}}

\DeclareMathOperator{\PM}{\mathrm{PM} } 

\def\.{\hskip.06cm}
\def\ts{\hskip.03cm}

\title{Normalized Matching Property for Competing Urn Models}

\date{\today}

 \author{Swee Hong Chan}
 \address[Swee Hong Chan]{Department of Mathematics,
 	Rutgers University, Piscataway, NJ, 08854.}
 \email{\texttt{sweehong.chan@rutgers.edu}}

\begin{document}

\begin{abstract}
We study the competing urn model in which \(m\) balls are placed independently
	into \(n\) urns according to (possibly distinct) ball distributions.
	Kahn and Neiman (2010) showed that, under identical ball distributions, the
	induced urn measure has \emph{conditional negative association} property and asked whether this
	remains true without assuming identical distributions.
We answer this in the affirmative by showing that the competing urn model
satisfies the \emph{normalized matching property}.
	This, in turn, implies conditional negative association for the induced urn measure  with non-identical
	ball distributions, resolving the question of Kahn and Neiman.
\end{abstract}


\maketitle

\section{Introduction} \label{s:intro}

In the classical \defnb{competing urns} experiment, $m$ balls are dropped 
independently into $n$ urns. 
The setup is simple to describe, but it already exhibits natural dependencies: 
more balls in one urn usually mean fewer elsewhere. 
While this negative correlation seems intuitively clear, proving that  
it continues to hold when  the state of some urns is revealed can be surprisingly 
delicate. 
This paper establishes this fact for generalized independent urn models.

Formally, we have a random function \. $\sigma:[m] \to [n]$ \. 
where \. $\sigma(1),\ldots, \sigma(m)$ \. are independent random variables. 
When the balls are identical (i.e. \. $\sigma(1),\ldots, \sigma(m)$ \.
are i.i.d),  we call this an  \defnb{ordinary urn model}.
When the balls are not identical, we call this  a \defnb{generalized independent urn model} instead (not to be confused with \emph{generalized P\'{o}lya urn model} in the literature, see e.g.~\cite{AK,HLS}).
Throughout this paper we will work with the generalized independent urn model, unless indicated otherwise.

\smallskip

Let   \. $B_j:=|\sigma^{-1}(j)|$ \.  be the number of balls in the $j$-th urn, and let $\mu$ be the law of \. $(B_1,\ldots, B_n)$\..
We call $\mu$ the (generalized independent) \defnb{urn enumeration measure}, which is a measure on $\nn^n$, where \ts $\nn := \{0,1,2,\ldots\}$ \ts is   the set of nonnegative integers.
 We denote by  \. $B_j^{\occ} \in \{0,1\}$  \.  the indicator for the occupation of the $j$-th urn, 
 and let  \. $\mu^{\occ}$ \. the law of \. $(B_1^{\occ},\ldots, B_n^{\occ})$\..
 We call \ts $\mu^{\occ}$ \ts the (generalized independent) \defnb{urn occupation measure}, which is  both a measure on $\nn^n$ and on $\{0,1\}^n$.
 In this paper, we investigate  negative dependence properties of these two urn measures, as detailed below.

\smallskip

\subsection{Conditional negative correlation and association}

We use boldface letters (e.g. $\bx$)  to denote  vectors in $\nn^n$,
and we write \ts $\bx=(\ax_1,\ldots, \ax_n)$, where \ts $\ax_j \in \nn$ \ts denotes the $j$-th entry of $\bx$.
We use $\leq $ to denote the  coordinatewise product order on $\nn^n$.
Let $\mu$ be a measure on $\nn^n$.
An event \ts $\cA \subseteq \nn^n$ \ts is \defnb{increasing} if,  for every \ts $\bx \in \cA$ \ts  and \ts $\by \in \nn^n $ \ts such that \. $\bx \leq \by$\., we have $\by \in \cA$.

For events $\cA,\cB \subseteq \nn^n$, we say that $\cA$ is \defnb{negatively correlated} with $\cB$, denoted by \ts $\cA \downarrow \cB$, if \. $\mu(\cA \cap \cB) \leq \mu(\cA) \mu(\cB)$.
Intuitively, this means that  conditioning on one  event does not increase the probability of the other event (throughout this paper we assume that any conditioning event we consider
has positive probability).
For  $\nn$-valued random variables $X,Y$, we write \. $X \downarrow Y$ \. if \. $ \{X \geq s\}  \. \downarrow  \. \{Y \geq t\}$ for all $s,t \in \nn$.
\defnb{Positive correlation} and  \. $X \uparrow Y$ \. are defined analogously.
We say that \ts $j \in [n]$ \ts  \defnb{affects} \ts $\cA \subseteq \nn^n$ \ts if there exist \ts $\bx \in \cA$ \ts and \ts $\by \in \nn^n \setminus \cA $ \ts such that 
\ts $\ax_i=\ay_i$ \ts for all \ts $i\in [n] \setminus \{j\}$\ts.
For events \ts $\cA,\cB \subseteq \nn^n$\ts, we write \ts $\cA \perp \cB$ \ts if no coordinate affects both $\cA$ and $\cB$. 
We say that a measure $\mu$ on $\nn^n$ is \defnb{negatively associated} if
\begin{equation}\label{eq:NA}
	\tag{NA}
 \text{\ts $\cA \downarrow \cB$ \ts whenever \ts $\cA,\cB \subseteq \nn^n$ \ts are increasing events and \ts $\cA \perp \cB$\ts.} 
\end{equation}
We say that $\mu$ is \defnb{negatively correlated}  if 
\begin{equation}\label{eq:NC}\tag{NC}
\text{$\ax_i \downarrow \ax_j$ \ 
whenever \ts $j \neq i$\ts, where $\bx \in \nn^n$ is randomly  sampled using $\mu$. }
\end{equation}
It follows from the definitions that \eqref{eq:NA} implies \eqref{eq:NC}.
On the other hand, \eqref{eq:NC} does not imply \eqref{eq:NA}, see e.g. \cite[Rem.~2.5]{JP}.

In this paper we study whether negative dependence is preserved under conditioning.
We say that $\mu$ is \defnb{conditionally negatively associated} \eqref{eq:CNA}  if
 any measure obtained from $\mu$ by conditioning the values of some of the variables is \eqref{eq:NA}.
That is, for all \ts $a_1,\ldots, a_n \in \nn$ \ts and $S \subseteq [n]$, 
\begin{equation}\label{eq:CNA} \tag{CNA}
\text{ the measure } \  \mu( \cdot \mid  \ax_i = a_i  \ \text{for all} \ i \in S)   \ \text{ satisfies } \ \eqref{eq:NA}.
\end{equation}
Similarly, $\mu$ is \defnb{conditionally negatively correlated} \eqref{eq:CNC} if, for all \ts $a_1,\ldots, a_n \in \nn$ \ts and $S \subseteq [n]$, 
\begin{equation}\label{eq:CNC} \tag{CNC}
	\text{ the measure } \  \mu( \cdot \mid  \ax_i = a_i  \ \text{for all} \ i \in S)   \ \text{ satisfies } \ \eqref{eq:NC}.
\end{equation}
It follows from the definitions that \eqref{eq:CNA} implies \eqref{eq:CNC}.
On the other hand, it is  an important  open problem posed by Pemantle  to determine if \eqref{eq:CNC} implies \eqref{eq:CNA}~\cite[Conj.~2]{Pem}.
We refer the reader to \cite{Pem} for an excellent survey of the background and motivations behind negative correlation.


\smallskip

\subsection{Back to urn measures}

The study of negative dependence in urn models was initiated 
 by Dubhashi and Ranjan~\cite{DR}, who showed that the  urn occupation 
measures satisfy \eqref{eq:NA}.
Building on this, and  motivated by the possibility that \eqref{eq:CNC} might fail to imply 
\eqref{eq:CNA}, Kahn and Neiman~\cite{KN} showed that both the generalized 
occupation and enumeration measures satisfy \eqref{eq:CNC}~\cite[Thm~14]{KN}. 
%
Thus, a negative answer to the following question of Kahn and Neiman would resolve  Pemantle’s problem.

%
%

\smallskip

\begin{question}[{\cite[Question~13]{KN}}]\label{quest:CNA}
	Are generalized independent urn measures \eqref{eq:CNA}?
\end{question}

\smallskip
Kahn and Neiman also proved that the \emph{ordinary} urn measures 
satisfy~\eqref{eq:CNA}~\cite[Thm.~1]{KN}, while noting that their method 
does not extend to the generalized setting. 
In this paper we answer this question affirmatively: both  urn occupation and enumeration measures satisfy \eqref{eq:CNA}.


\smallskip

\begin{thm}\label{thm:CNA-occ}
The generalized independent urn occupation measure 	$\mu^{\occ}$ satisfies \eqref{eq:CNA}.
\end{thm}

\smallskip

\begin{thm}\label{thm:CNA-ball}
	The generalized independent urn enumeration measure 	$\mu$ satisfies \eqref{eq:CNA}.
\end{thm}

\smallskip

On the one hand, these theorems resolve the question of Kahn and Neiman, which they noted
as the most interesting question left open by their work~\cite[\S5]{KN}.
On the other hand, this also means that the generalized independent urn measures cannot be used to settle  the open problem of Pemantle.

It is perhaps worth noting that generalized independent urn models are known to \emph{not} satisfy 
ultra-log-concavity or the Rayleigh property (see \S\ref{ss:lc} and \S\ref{ss:Rayleigh}), two 
properties closely related to \eqref{eq:CNA}. 
This might suggest that the theorems proved here are subtler than they 
first appear.

We prove Theorem~\ref{thm:CNA-occ} in \S\ref{subsec:main-thm-occ}, and Theorem~\ref{thm:CNA-ball} in \S\ref{subsec:main-thm-ball}.
The key step in both cases is to show that  generalized independent urn models satisfy the \defnb{normalized matching property}, a structural condition that, together with \eqref{eq:CNC}, is known to imply \eqref{eq:CNA}.

\smallskip

\subsection{Normalized matching property}\label{ss:NMP}

Let $\nu$ 
be a measure on  $[m]$, and 
let \ts $Z \subseteq [m]$ \ts  be sampled according to $\nu$.
In the context of this paper, we will think of  $Z$ as the set of balls that end up in the $n$-th urn.
The measure  $\nu$ satisfies the \defnb{normalized matching property} if
\begin{equation}\label{eq:NMP}\tag{NMP}
	\nu(\cA \mid  |Z|=k) \ \leq \  \nu(\cA \mid  |Z|=k+1) \quad \text{ for all $k\geq 1$ and increasing  $\cA \subseteq 2^{[m]}$}.
\end{equation}
That is, the measure \. $\nu(\cdot \mid |Z|=k+1)$ \. \emph{stochastically dominates} \. $\nu(\cdot \mid |Z|=k)$\..
This terminology comes from the corresponding notion in graph theory, and the relationship between the two will be discussed in greater detail in  \S\ref{ss:bipartite}. 

In this paper, we will prove \eqref{eq:NMP} for the following two measures.
Let $d \in \{0,1,\ldots, n-1\}$, and let  \ts $Q_d^{\occ}$ \ts denote the event that the first $d$ urns are nonempty,   
\begin{equation*}
	Q_d^{\occ} \ := \ \{ B_1 \geq 1,\ldots, B_d\geq 1  \} \ = \ \{ B_1^{\occ} = 1,\ldots, B_d^{\occ} = 1  \}.
\end{equation*}
We denote by \. $\nu_d^{\occ}$ \.  the distribution of $\sigma^{-1}(n)$  conditioned on  $Q_d^{\occ}$.  
That is to say, \. $\nu_d^{\occ}$ \. is the law of the set of balls in the $n$-th given that the first $d$ urns are non-empty.

\smallskip

\begin{thm}\label{thm:NMP-occ}
	The measure $\nu_d^{\occ}$ satisfies \eqref{eq:NMP} \. for all \. $d\in \{0,1,\ldots,n-1\}$\..
\end{thm}

\smallskip

For the second measure, let  \. $\ba=(\aa_1,\ldots, \aa_d) \in \nn^d$, and let  $Q_{\ba}$ be  
\[ Q_{\ba} \ := \ \{ B_1 =\aa_1,  \. \ldots, \. B_d=\aa_d  \},\]
the event that the number of balls in the first $d$ urns are specified by $\ba$.
We denote by \. $\nu_{\ba}$ \. the distribution  of $\sigma^{-1}(n)$  conditioned on the event $Q_{\ba}$.  

 \smallskip

\begin{thm}\label{thm:NMP-ball}
	The measure $\nu_{\ba}$ satisfies \eqref{eq:NMP} for all \. $\ba \in \nn^d$.
\end{thm}

\smallskip

Whether the  measures \ts  $\nu_d^{\occ}$ \ts and \ts $\nu_{\ba}$ \ts satisfy \eqref{eq:NMP} was first asked by Kahn--Neiman~\cite[Question~17]{KN}. 
They noted, though without writing out the details, that this fact already suffices to deduce Theorems~\ref{thm:CNA-occ} and \ref{thm:CNA-ball}.
 We will prove  Theorem~\ref{thm:NMP-occ} in \S\ref{ss:proof-NMP-occ} and Theorem~\ref{thm:NMP-ball} in \S\ref{ss:proof-NMP-ball}, and then present in \S\ref{sec:FM-argument} the 
 omitted argument that leads from these results to 
 Theorems~\ref{thm:CNA-occ} and~\ref{thm:CNA-ball}.



\medskip

The structure of the paper is as follows.  
In \S\ref{s:prelim} we review preliminaries on normalized matching 
properties and graph orientations.  
In \S\ref{sec:proof-orient-occ} we prove the main technical lemma (Lemma~\ref{lem:orient-occ}) that will be  used in the proof of Theorem~\ref{thm:NMP-occ}.  
In \S\ref{sec:proof-orient-ball} we prove the main technical lemma (Lemma~\ref{lem:orient-ball}) that will be  used in the proof of Theorem~\ref{thm:NMP-ball}.  
In \S\ref{s:proof-H-NMP} we give the proofs of Theorems~\ref{thm:NMP-occ} 
and~\ref{thm:NMP-ball}.  
In \S\ref{sec:FM-argument} we prove Theorems~\ref{thm:CNA-occ} 
and~\ref{thm:CNA-ball}.  
Finally, in \S\ref{s:finrem} we conclude with final  remarks and open problems.

\medskip

\section{Preliminaries and Setup}\label{s:prelim}


\subsection{Normalized matching property for bipartite graphs}\label{ss:bipartite}
In this section, we discuss the normalized matching property for  bipartite graphs and show how it can be applied to derive \eqref{eq:NMP} in the context of  urn measures.

Let \. $G=(V_1, V_2, E)$ \. be  a  bipartite graph
with bipartition $V_1$ and $V_2$, and edge set $E$.
Let \. $f:V_1 \cup V_2  \to \Rb_{\geq 0}$ 
be a nonnegative weight function on the vertices of $G$.
For any  \ts $U\subseteq V_1 \ts \cup  \ts V_2$, we write \. $f(U):=\sum_{v \in U} f(v)$. 
We assume that  $f$ is \defnb{balanced}, i.e. \. $f(V_1)=f(V_2)$.
We  denote by \. $N(U)$ \. the set of vertices  adjacent to $U$.
The following are equivalent formulations of the \defnb{normalized matching property for bipartite graphs} with balanced weight functions:
\begin{enumerate}
	\item \label{it:a} For all  \. $U \subseteq V_1$\.,   
	\begin{equation*}\label{eq:HMC}\tag{HMC}
	   {f(U)}  \  \leq \  {f(N(U))}. 
	\end{equation*}
	\item \label{it:b} For all independent set $U \subseteq V_1 \cup V_2$,
	\[   f(U\cap V_1) \. + \. f(U \cap V_2) \  \leq \ f(V_1). \]
	\item \label{it:c} There exists \. $\omega:E \to \Rb_{\geq 0}$ \. such that, for all $v \in V_1 \cup V_2$,
	\begin{equation*}\label{eq:NMP-b}\tag{NMP-B}
		\sum_{e \in E \text{ incident to $v$}} \omega(e) \ = \  f(v). \
	\end{equation*}
\end{enumerate}
The normalized matching property was first introduced by Graham and Harper~\cite{GH} in the form of condition~\eqref{it:a}, as a weighted generalization of Hall’s classical matching theorem~\cite{Hall}.
Condition~\eqref{it:b} is a special case of the classical LYM property (also known as the YBLM property)~\cite{Yam,Bol,Lub,Mes}; see \eqref{eq:LYM} in \S\ref{ss:base} for its formulation in ranked posets.
In what follows, we will primarily use condition~\eqref{it:c}, which is the most convenient for our purposes.
A proof of the equivalence of these three conditions via the max–flow min–cut theorem is given in~\cite{Kle}.\footnote{Although the results in~\cite{Kle} are stated for the uniform weight function, the arguments extend without difficulty to the general case considered here.}
Note that if two nonnegative functions $f_1$ and $f_2$ satisfy \eqref{eq:NMP-b}, then so does their sum $f_1+f_2$. We will make repeated use of this property throughout the paper.


\smallskip

We now construct a weighted bipartite graph arising from the urn model in \S\ref{s:intro}.
Recall that $[m]$ is the set of balls and $[n]$ is the set of urns in the generalized independent urn model.

\smallskip
\begin{definition}[Urn incidence graph]\label{def:H}
Let \. $X\subseteq [m]$ \. be an  odd set, and set \. $k:=\tfrac{|X|-1}{2}$.
Define the bipartite graph \. $H_X:=(W_1,W_2,F):=(W_1(X), W_2(X), F(X))$ \.  by
\begin{align*}
	& W_{1} \ := \ \binom{X}{k}, \quad W_{2} \ := \  \binom{X}{k+1}, \\
	& F \ := \ \{ \. (S,T) \in W_{1} \times W_{2}  \. \mid \. S \subset T \. \}.
\end{align*}
Let 
$\nu$ be a probability measure on $[m]$ (corresponding to the random balls ending in the
$n$-th urn).
We define a weight function  \. $g_{\nu}:=g_{\nu,X}$ \. on the vertices of  $H_X$  by
\[  g_{\nu}(S) \ := \ \nu(S) \. \nu(X\setminus S).  \]
By construction, $g_{\nu,X}$ is a balanced weight function.
\end{definition}

\smallskip

The next two  propositions show that the weighted bipartite graphs associated with the measures $\nu_d^{\occ}$  and $\nu_{\ba}$  from \S\ref{s:intro} satisfy \eqref{eq:NMP-b}.

\smallskip

\begin{prop}\label{prop:H-occ}
		For any odd set \ts $X \subseteq [m]$ \ts and   \ts $d \in \{0,1,\ldots,n-1\}$\ts,
		the 
  bipartite graph \. $H_X$ \. with the weight function  \. $g_{\nu_d^{\occ}}$ \. satisfies the normalized matching property \eqref{eq:NMP-b}.
\end{prop}

\smallskip

\begin{prop}\label{prop:H-ball}
		For any odd set \ts $X \subseteq [m]$ \ts and   \ts $\ba \in \nn^d$\ts,
the 
bipartite graph \. $H_X$ \. with the weight function  \. $g_{\nu_{\ba}}$ \. satisfies the normalized matching property \eqref{eq:NMP-b}.
\end{prop}

\smallskip


We will prove Proposition~\ref{prop:H-occ} in \S\ref{ss:proof-H-occ} and Proposition~\ref{prop:H-ball} in \S\ref{ss:proof-H-ball}, with the next two sections devoted to laying the groundwork. In particular, the proofs make use of graph orientations, reviewed in the next subsection.


\medskip

\subsection{Admissible graph orientations}\label{ss:orient}

Let \. $X \subseteq [m]$ \. be an odd set. 
Let $G$ be an undirected graph (possibly with parallel edges and loops) with $n$ vertices $v_1,\ldots, v_n$
and with $m$ edges $e_1,\ldots, e_m$.
Throughout the paper we assume that $G$ satisfies
\begin{equation}\label{eq:X-G}
X \. = \. \{ \. j \in [m] \. \mid \. e_j \ \text{ incident to } \. v_n  \. \text{ in $G$}\}. 
\end{equation}
That is to say, $X$ encodes the edge set at the vertex with the largest label.

An \defnb{orientation} \(O\) of \(G\) is an assignment of a direction to each edge of \(G\). 
For \(i \in [n]\), let \(\Out_O(v_i)\) and \(\In_O(v_i)\) denote the sets of outgoing and incoming 
edges of \(v_i\) in \(O\), respectively, and write 
\. $\outdeg_O(v_i) := |\Out_O(v_i)|$ \. and 
\. 
$\indeg_O(v_i) := |\In_O(v_i)|$\..

\smallskip

\begin{definition}[$d$-admissible orientations]
For any nonnegative integer $d$, 
an orientation $O$ of $G$ is \defnb{$d$-admissible} if 
\[ \text{$\outdeg_O(v_i) \geq 1$ \. and \.  \. $\indeg_O(v_i) \geq 1$ \quad   for all $i \in [d]$.} \]
For $S \subseteq X$, define
\begin{align*}
	\aM_{G,d}^{\occ}(S) \ := \ 
	\begin{array}{l}
	\text{number of $d$-admissible orientations $O$ of $G$ with  }\\[4pt]
	\text{\. $\Out_O(v_n)=\{e_j \. \mid \. j \in S\}$ \ and \  $\In_O(v_n)=\{e_j \. \mid \. j \in X\setminus S$\}.} 
	\end{array}
\end{align*}
Here the  notation 
$\occ$ anticipates its connection to the urn occupancy measure, which will be made explicit in \S\ref{ss:proof-H-occ}.
By construction,  \ts $\aM_{G,d}^{\occ}$ \ts is a balanced weight function for vertices of $H_X$.
\end{definition}

\smallskip

\begin{definition}[$\ba$-admissible orientations]
For any  \. $\ba=(\aa_1,\ldots, \aa_d) \in \nn^d$\., an 
 orientation $O$ of $G$ is 
\defnb{$\ba$-admissible} if 
\[ \outdeg_O(v_i)=\indeg_O(v_i)=\aa_i \quad   \text{ for all }  i \in [d].\]
For $S \subseteq X$, define 
\begin{align*}
	\aM_{G,\ba}(S) \ := \ 
		\begin{array}{l}
		\text{number of $\ba$-admissible orientations $O$ of $G$ with  }\\[4pt]
		\text{\. $\Out_O(v_n)=\{e_j \. \mid \. j \in S\}$ \ and \  $\In_O(v_n)=\{e_j \. \mid \. j \in X\setminus S$\}.} 
	\end{array} 
\end{align*}
By construction, \ts $\aM_{G,\ba}$ \ts is a balanced weight function for vertices of $H_X$.
\end{definition}

The connections between admissible orientations and the normalized matching 
property will be established in \S\ref{sec:proof-orient-occ} and \S\ref{sec:proof-orient-ball}.


\medskip

\section{Normalized matching property and $d$-admissible orientations}\label{sec:proof-orient-occ}

This section is devoted to proving the following lemma, which forms the 
basis for Proposition~\ref{prop:H-occ}, and hence Theorems~\ref{thm:CNA-occ} 
and~\ref{thm:NMP-occ}.
Throughout this section, \. $X \subseteq [m]$ \. is an odd set,
\. $G$  \. is a graph with vertices $v_1,\ldots, v_n$ and edges $e_1,\ldots, e_m$ satisfying  \eqref{eq:X-G}, and  \. $d \in \{0,1\ldots,n-1\}$\..
Recall the definitions of \ts $H_X$ \ts  and \ts $\aM_{G,d}^{\occ}$ \ts  from \S\ref{s:prelim}.

\smallskip

\begin{lemma}\label{lem:orient-occ}
The   bipartite graph \. $H_{X}$ \. with the weight function \. $\aM_{G,d}^{\occ}$ \. satisfies the normalized matching property \eqref{eq:NMP-b}.
\end{lemma}

\smallskip

We now build toward the proof of Lemma~\ref{lem:orient-occ}.


\subsection{Deletion and contraction recurrences}\label{subsec:rec-loop}
In this section we collect various  recursion formulas for the function \.$\aM_{G,d}^{\occ}$ \. that will be used in the proof of Lemma~\ref{lem:orient-occ}.
Let $e$ be an edge of $G$.  
The \defnb{edge deletion} $G-e$ is the graph obtained by removing $e$,  
and the \defnb{edge contraction} $G/e$ is the graph obtained by identifying the endpoints of $e$ and then deleting $e$.  
Since the definition of \. $\aM_{G,d}^{\occ}(S)$ \. depends on the ordering of vertices 
$v_1,\ldots,v_n$ and edges $e_1,\ldots,e_m$ of $G$, we must prescribe compatible orderings 
for $G-e$ and $G/e$.  
For $G-e$, we adopt the ordering inherited from $G$, with the convention that the deleted edge is always the last edge $e_m$.  
For $G/e$, we retain the vertex and edge ordering from $G$, with the additional rule that the contracted vertex inherits the label of the smaller endpoint. In other words, if $e$ joins $v_i$ and $v_j$ with $i \leq j$, then the vertices of $G/e$ are $v_1', \ldots, v_{n-1}'$, where $v_k' = v_k$ for $k < j$, and $v_k' = v_{k-1}$ for $k \geq j$. In the resulting graph $G/e$, the set $X$ is understood to be the set of indices of edges incident to the last vertex, namely $v_{n-1}' = v_n$.

\smallskip

\begin{lemma}[Loop-edge recurrences]\label{lem:N-occ-loop}
Let $e=e_m$ be a loop edge of $G$.  
	Then, for any  $S \subseteq X$:
	\begin{enumerate}
		\item \label{it:N-occ-loop-1} If  the loop $e$ is at  $v_d$, then 
		\[ \aM_{G,d}^{\occ}(S) \ = \ \aM_{G-e,d-1}^{\occ}(S).  \]
		\item \label{it:N-occ-loop-2} If  the loop $e$ is at $v_{d+1}$, then 
\[ \aM_{G,d}^{\occ}(S) \ = \  \aM_{G-e,d}^{\occ}(S).  \]
		\item \label{it:N-occ-loop-3} If  the loop $e$ is at  $v_{n}$, then 
\[ \aM_{G,d}^{\occ}(S) \ = \ 0.  \]
	\end{enumerate} 
\end{lemma}

\smallskip

\begin{proof}
For \eqref{it:N-occ-loop-1}, if $e$ is a loop at $v_d$, then every orientation of $e$ ensures  
$\indeg_O(v_d)\geq 1$ and $\outdeg_O(v_d)\geq 1$, so the recursion holds.  
For \eqref{it:N-occ-loop-2}, if $e$ is a loop at $v_{d+1}$, then the orientation of $e$ does not affect 
any of the $d$-admissibility conditions, hence the recursion.  
For \eqref{it:N-occ-loop-3}, if $e$ is a loop at $v_n$, then any orientation of $e$ makes it impossible to satisfy simultaneously  
$\Out(v_n)=\{e_j \mid j \in S\}$ and $\In(v_n)=\{e_j \mid j \in X\setminus S\}$.  
Thus $\aM_{G,d}^{\occ}(S)=0$.
\end{proof}

\smallskip

For an edge \. $e = \{v_i, v_j\}$ \. of $G$, we say that $e$ is of \defnb{type~A} if \. $i,j \leq d$\., of \defnb{Type~B} if \. $i \leq d < j < n$\., and of \defnb{Type~C} if \. $d < i,j < n$\..

\smallskip

\begin{lemma}[Type A recurrences]\label{lem:N-occ-1}
Let $e = e_m$ be an edge of $G$, 
and suppose that $e$ is incident to $v_{d-1}$ and $v_d$. 
Then, for any  $S \subseteq X$, 
\[
\aM_{G,d}^{\occ}(S) \;=\;
\begin{cases}
	\aM_{G-e,d}^{\occ}(S) \;+\; \aM_{G/e,d-1}^{\occ}(S), 
	& \text{if $\deg_G(v_{d-1}) \geq 2$ and $\deg_G(v_d) \geq 2$}, \\[6pt]
	0, & \text{otherwise}.
\end{cases}
\]
\end{lemma}

\smallskip

\begin{proof}
If either $\deg_G(v_{d-1})=1$ or $\deg_G(v_d)=1$, then in any orientation one of indegree/outdegree at that vertex is zero, so no $d$-admissible orientations exist.
Thus, we may assume that $\deg_G(v_{d-1}) \geq 2$ and $\deg_G(v_d) \geq 2$.

	The argument now follows the standard deletion–contraction approach.	
	We  partition the orientations counted by $\aM_{G,d}^{\occ}(S)$ into three  families of orientations $O$:
	\begin{enumerate}
		\item\label{it:N-occ-V1-1}  $O-e$ \ts (the orientation obtained by removing the edge $e$ from the digraph) is $d$-admissible, and  $e$ is oriented from $v_{d-1}$ to $v_d$ in $O$.
		
		\item\label{it:N-occ-V1-2}   	  $O-e$ \ts is $d$-admissible, and  $e$ is oriented from $v_{d}$ to $v_{d-1}$ in $O$.
		
				\item \label{it:N-occ-V1-3}  	  $O-e$ \ts is not $d$-admissible.
In this case, one of the following holds:
\begin{enumerate}[label=(\alph*)]
	\item \label{it:A-a} $\outdeg_{O-e}(v_{d-1}) = 0$;
	\item \label{it:A-b} $\indeg_{O-e}(v_{d-1}) = 0$;
	\item \label{it:A-c} $\outdeg_{O-e}(v_d) = 0$; or
	\item \label{it:A-d} $\indeg_{O-e}(v_d) = 0$.
\end{enumerate}
  These subcases are not mutually exclusive, but in each the orientation of $e$ is forced: 
$e$ is oriented from $v_{d-1}$ to $v_d$ in \ref{it:A-a} and \ref{it:A-d}, 
and from $v_d$ to $v_{d-1}$ in \ref{it:A-b} and \ref{it:A-c}.
	\end{enumerate}

	\smallskip
	
	\noindent\textbf{Deletion part.}  
	Orientations in \eqref{it:N-occ-V1-1} are in bijection with those counted by 
	$\aM_{G-e,d}^{\occ}(S)$, obtained by deleting $e$. 
	
	\smallskip

	\noindent\textbf{Contraction part.}  
	Orientations in \eqref{it:N-occ-V1-2} and \eqref{it:N-occ-V1-3} are in bijection with those counted by 
	$\aM_{G/e,d-1}^{\occ}(S)$, obtained by contracting $e$. 
	To see this, note that every orientation counted by 
	$\aM_{G/e,d-1}^{\occ}(S)$ belongs to one of the following families:
	\begin{enumerate}[label=(\arabic*),start=0]
		\item\label{it:A-0} The edges originally incident to $v_{d-1}$ include both an incoming and an outgoing edge at $v_{d-1}$, and similarly for $v_d$.
		\item\label{it:A-1} All edges originally incident to $v_{d-1}$ are oriented into $v_{d-1}$.
		\item\label{it:A-2} All edges originally incident to $v_{d-1}$ are oriented out of $v_{d-1}$.
		\item\label{it:A-3} All edges originally incident to $v_d$ are oriented into $v_d$.
		\item\label{it:A-4} All edges originally incident to $v_d$ are oriented out of $v_d$.
	\end{enumerate}
	The cases \ref{it:A-1}--\ref{it:A-4} are not mutually exclusive, and the edge sets described there are 
	nonempty since $\deg_G(v_{d-1}), \deg_G(v_d)\geq 2$. 
In each such case, one can recover a $d$-admissible orientation of $G$ by adding back the edge $e$, with $e$ directed from $v_{d-1}$ to $v_d$ in \ref{it:A-a} and \ref{it:A-d}, 
and from $v_d$ to $v_{d-1}$ in \ref{it:A-b} and \ref{it:A-c}. 
	Now orientations in \eqref{it:N-occ-V1-2} correspond exactly to family \ref{it:A-0}, 
	while orientations in \eqref{it:N-occ-V1-3} correspond to families \ref{it:A-1}--\ref{it:A-4}. 
	
	\smallskip

	Together these cases exhaust all orientations,  and the proof is complete.
\end{proof}

\smallskip

\begin{lemma}[Type B recurrences]\label{lem:N-occ-2}
Let $e = e_m$ be an edge of $G$, and let $1 \leq d < n-1$. 
Suppose that $e$ is incident to $v_{d}$ and $v_{d+1}$. 
Then, for any $S \subseteq X$,  
	\[   \aM_{G,d}^{\occ}(S) \ = \ 
\begin{cases}
\aM_{G-e,d}^{\occ}(S)  \. + \. \aM_{G-e,d-1}^{\occ}(S) &\text{ if } \deg_G(v_d)\geq 2,\\
0 & \text{ if } \deg_G(v_d)=1.
\end{cases} \]
\end{lemma}

\smallskip

\begin{proof}
	If  $\deg_G(v_d)=1$, then in any orientation either the indegree or the outdegree of $v_d$ is zero, so no $d$-admissible orientations exist. Thus we may assume that  \. $\deg_G(v_d)\geq 2$\., in which case there exists another edge $e'$ that is incident to $v_d$.
We now partition the orientations counted by $\aM_{G,d}^{\occ}(S)$ into two families:
\begin{enumerate}
	\item Orientations $O$ in which $e$ is oriented out of $v_d$ and $e'$ is also oriented out of $v_d$, 
	or $e$ is oriented into $v_d$ and $e'$ is also oriented into $v_d$.  
	These orientations are in bijection with those counted by $\aM_{G-e,d}^{\occ}(S)$, 
	obtained by deleting the edge $e$. 
	
	\item Orientations $O$ in which $e$ is oriented out of $v_d$ while $e'$ is oriented into $v_d$, 
	or $e$ is oriented into $v_d$ while $e'$ is oriented out of $v_d$.  
	These orientations are in bijection with those counted by $\aM_{G-e,d-1}^{\occ}(S)$, 
	again obtained by deleting the edge $e$. 
\end{enumerate}
This completes the proof.
\end{proof}

\smallskip

\begin{lemma}[Type C recurrences]\label{lem:N-occ-3}
	Let $e = e_m$ be an edge of $G$, and let  $d < n-2$. 
	Suppose that $e$ is incident to $v_{d+1}$ and $v_{d+2}$. 
	Then, for any $S \subseteq X$,  	
\[   \aM_{G,d}^{\occ}(S) \ = \ 2 \aM_{G-e,d}^{\occ}(S). \]
\end{lemma}

\smallskip

\begin{proof}
	Either orientation of the edge $e$ leaves the $d$-admissibility of the rest of the orientation unchanged. 
	Thus each orientation of $G-e$ extends in exactly two ways to an orientation of $G$, 
	and the result follows.
\end{proof}

\smallskip

\subsection{Base case of Lemma~\ref{lem:orient-occ}}\label{ss:base} 

We now present a special case of Lemma~\ref{lem:orient-occ}, which will serve as the foundation for our inductive proof of the full lemma in the next subsection.

\begin{lemma}\label{lem:N-occ-base}
	Suppose that every edge of $G$ is incident to $v_n$.
	Then the bipartite graph $H_X$ with the weight function $\aM_{G,d}^{\occ}$ satisfies the normalized matching property \eqref{eq:NMP-b}.
%
\end{lemma}

\smallskip 

\begin{proof}
	First note that if $G$ contains a loop, then by Lemma~\ref{lem:N-occ-loop}\eqref{it:N-occ-loop-3} 
	we have $\aM_{G,d}^{\occ}\equiv 0$, and the claim follows trivially. 
	Thus we may assume that $G$ has no loops.
	
	For $1 \leq i \leq n$, let $X_i \subseteq [m]$ be the set of indices of edges incident to both $v_i$ and $v_n$. 
	Set $k := \tfrac{|X|-1}{2}$. 	
By definition, for any subset $S \subseteq X$ of size $k$ or $k+1$, we have
$\aM_{G,d}^{\occ}(S)$ is the indicator that $S$ intersects each part $X_i$ nontrivially but not completely, i.e.
	\begin{align*}
		 \aM_{G,d}^{\occ}(S) \ = \ 
		 \begin{cases}
 1  & \text{ if } \  0< |S \cap X_i|  < |X_i| \. \text{ for } 1 \leq i \leq d,\\
0 & \text{ otherwise}.
		 \end{cases}
		\end{align*}

	Let \. $H':=(W_1',W_2',E')$ \.  be the subgraph of $H_X$
	induced by the vertices  $S$ with \. $\aM_{G,d}^{\occ}(S)\geq 1$\., i.e.
	\begin{align*}
W_1' \ &:= \ \bigg\{ \. S \in \binom{X}{k} \. \mid \.   0< |S \cap X_i|  < |X_i| \. \text{ for } 1 \leq i \leq d \. \bigg \},\\
W_2' \ &:= \ \bigg\{ \. S \in \binom{X}{k+1} \. \mid \.   0< |S \cap X_i|  < |X_i| \. \text{ for } 1 \leq i \leq d \. \bigg \},\\
E' \ &:= \  \{ \. \{S,T\} \in W_{1}' \times W_{2}'  \. \mid \. S \subset T \. \}.
	\end{align*}
By the equivalent formulations of the normalized matching property for weighted bipartite graphs 
(see \S\ref{ss:bipartite}), it suffices to verify that $H'$ satisfies the Hall's  condition \eqref{eq:HMC}:
	\begin{align*}
		|U| \ \leq \ |N(U)| \quad \text{for all } U \subseteq W_1'.
	\end{align*}
We will establish this as a special case of a classical result of Griggs~\cite{Gri}. 
	A ranked poset $P$ is said to have  the \defnb{LYM property} (also known as \defnb{YBLM property}) if every antichain $U$ of $P$ satisfies 
	\begin{equation}\label{eq:LYM}\tag{LYM}
	 \sum_{x \in U} \frac{1}{w(x)} \ \leq \ 1,  
	\end{equation}
where $w(x)$ is the number of elements in  $P$ having the same rank as $x$.

\smallskip

\begin{thm}[{\cite[Thm.~3.2]{Gri}}]\label{thm:Griggs}
	Let $X$ be a finite set partitioned into $n$ parts $X_1,\ldots,X_n$.  
	For each $i \in [n]$, let \. $I_i \subseteq \{0,1,\ldots,|X_i|\}$ \. be an arithmetic progression. 
	Define $P$ to be the poset with ground set
	\[
	\bigl\{\, S \subseteq X \, \mid \,  |S \cap X_i| \in I_i \, \text{ for } 1 \leq i \leq n \,\bigr\},
	\]
	ordered by inclusion.  
	Then $P$ satisfies \eqref{eq:LYM}.
\end{thm}


\smallskip

We now apply Theorem~\ref{thm:Griggs} with \. $I_i$ \. ($i \in [n]$) \. given by
\begin{align*}
	I_i \  := \  \begin{cases}
		\{1,\ldots, |X_i|-1\} & \text{ for } \ 1 \leq i \leq d,\\
				\{0,1,\ldots, |X_i|\} & \text{ for } \ i \geq d+1.
	\end{cases}
\end{align*}
It  follows that $H'$ is precisely the Hasse diagram of $P$ restricted to the  elements of rank $k$ or $k+1$.
Now, for any subset $U \subseteq W_1'$, it follows from applying \eqref{eq:LYM} to \. $U \cup (W_2' \setminus N(U))$ \. that 
\[    \frac{|U|}{|W_1'|} \. + \. \frac{|W_2'|-|N(U)|}{|W_2'|}  \ \leq \ 1, \]
which is equivalent to \. $|U|\leq |N(U)|$\.,  since $|W_1'|=|W_2'|$ by construction. 
This establishes \eqref{eq:HMC} and completes the proof.
\end{proof}

\smallskip

\subsection{Proof of Lemma~\ref{lem:orient-occ}}
Let $G$ be the graph in the statement of the lemma. 
We proceed by induction on the number of edges $m$. 
Let $e$ be an arbitrary edge of $G$. By relabeling if necessary, 
we may assume $e = e_m$ is the last edge of $G$.
We distinguish cases according to the endpoints of  $e$.

%
%
%
%

\smallskip

\noindent
\textbf{Case 1: $e$ is a loop.}  
By Lemma~\ref{lem:N-occ-loop} the problem reduces to the case of $G-e$, 
which follows by the induction hypothesis. 
Hence we may assume that $e$ is a proper edge.

\smallskip

\noindent
\textbf{Case 2: both endpoints of $e$ lie in $\{v_1,\ldots,v_d\}$.}  
We may assume without loss of generality that $e$ is incident to 
$v_{d-1}$ and $v_d$. In this case, Lemma~\ref{lem:N-occ-1} reduces the 
problem to the cases of $G-e$ and $G/e$, which are handled by the 
induction hypothesis, together with the fact that the sum of nonnegative 
functions satisfying \eqref{eq:NMP-b} also satisfies \eqref{eq:NMP-b}.

\smallskip

\noindent
\textbf{Case 3: $e$ has one endpoint in $\{v_1,\ldots,v_d\}$ and the 
	other in $\{v_{d+1},\ldots,v_{n-1}\}$.}  
We may assume without loss of generality that $e$ is incident to $v_d$ 
and $v_{d+1}$. Then by Lemma~\ref{lem:N-occ-2} the problem reduces to the 
case of $G-e$, which is covered by the induction hypothesis, again using 
the closure of \eqref{eq:NMP-b} under addition.

\smallskip

\noindent
\textbf{Case 4: both endpoints of $e$ lie in $\{v_{d+1},\ldots,v_{n-1}\}$.}  
We may assume without loss of generality that $e$ is incident to 
$v_{d+1}$ and $v_{d+2}$. Then by Lemma~\ref{lem:N-occ-3} the problem 
reduces to the case of $G-e$, which follows by induction.

\smallskip

\noindent
\textbf{Case 5: every edge of $G$ is incident to $v_n$.}  
In this final case, the claim follows directly from 
Lemma~\ref{lem:N-occ-base}.  

\smallskip

The proof is now complete. \qed

\medskip

\section{Normalized matching property and $\ba$-admissible orientations}\label{sec:proof-orient-ball}

This section is devoted to proving the following lemma, which forms the 
basis for Proposition~\ref{prop:H-ball}, and hence Theorems~\ref{thm:CNA-ball} 
and~\ref{thm:NMP-ball}.
Throughout this section, \. $X \subseteq [m]$ \. is an odd set,
\. $G$  \. is a graph with vertices $v_1,\ldots, v_n$ and edges $e_1,\ldots, e_m$ satisfying  \eqref{eq:X-G}. We fix    \. $d \in \{0,1\ldots,n-1\}$ \.  and $\ba=(\aa_1,\ldots, \aa_d) \in \nn^d$.
Recall the definitions of \ts $H_X$ \ts  and \ts $\aM_{G,\ba}$ \ts  from \S\ref{s:prelim}.

\smallskip

\begin{lemma}\label{lem:orient-ball}
The   bipartite graph \. $H_{X}$ \. with the weight function \. $\aM_{G,\ba}$ \. satisfies the normalized matching property \eqref{eq:NMP-b}.
\end{lemma}

\smallskip

Throughout this section we assume that $\deg_G(v_i) = 2\aa_i$ for all 
$i \in [d]$; otherwise $\aM_{G,\ba} \equiv 0$, and Lemma~\ref{lem:orient-ball} 
holds trivially.
To prove Lemma~\ref{lem:orient-ball}, we begin with several recurrence 
identities for $\aM_{G,\ba}$.

\smallskip

\subsection{Recurrence identities}
Recall the definition of edge deletion $G-e$ from \S\ref{subsec:rec-loop}. The first recurrence identity concerns loop edges.
\smallskip

\begin{lemma}[Loop-edge recurrences]\label{lem:N-ball-loop}
Let $e = e_m$ be a loop at $v_d$ in $G$. 
	Then, for any  $S \subseteq X$:	
\[ \aM_{G,\ba}(S) \ = \   \aM_{G-e,\ba'}(S),  \quad \text{ where } \quad  \ba' \ := \ (\aa_1,\ldots, \aa_{d-1}, \aa_{d}-1) \in \nn^{d}.\] 

\end{lemma}

\smallskip

\begin{proof}
This follows immediately from the fact that the orientation of the loop $e$ 
contributes one incoming edge and one outgoing edge to $v_d$.
\end{proof}

\smallskip

The second recurrence identity requires some notation. 
Let $Y_d$ denotes the set of edges incident to $v_d$, and assume that 
no loops are incident to $v_d$. 
Recall from the assumption at the beginning of this subsection that  \. $|Y_d| = \deg_G(v_d) = 2\aa_d$ \. is even. 
We write $\PM(Y_d)$ for the set of \defnb{perfect matchings} of $Y_d$, 
that is, the set of fixed-point-free involutions on $Y_d$. 
Equivalently, each $\pi \in \PM(Y_d)$ pairs the edges of $Y_d$ into disjoint 2-element subsets.



\begin{figure}[ht!]
	\label{fig:split}
	\centering
	\includegraphics[scale=1]{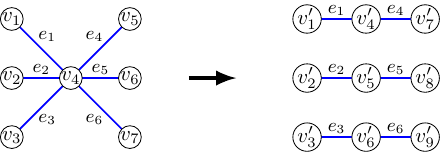}
	\caption{(Left) The graph $G$. (Right) The $\pi$-split of $v_d$ with $d=4$ and $\pi=(14)(25)(37)$.}
\end{figure}

For $\pi \in \PM(Y_d)$, we define $G[\pi]$ as the graph obtained from $G$ by 
\defnb{splitting} the vertex $v_d$ according to $\pi$. 
That is, delete $v_d$ and replace it with $\aa_d$ new vertices 
$v_d',\ldots,v_{d+\aa_d-1}'$, and for each pair \. $\{e,e'\} \in \pi$ \. 
reattach both $e$ and $e'$ to a distinct new vertex. 
Thus each new vertex becomes the common endpoint of exactly one pair. 
We will refer to this construction as the \defnb{$\pi$-split of $v_d$}; 
see Figure~\ref{fig:split} for an illustration.

Since the definition of $\aM_{G,\ba}$ depends on the ordering of vertices and edges, 
we extend the vertex order by inserting $v_d',\ldots,v_{d+\aa_d-1}'$ immediately 
after $v_{d-1}$, and relabel so that $v_i' = v_i$ for $i<d$ and 
$v_{i+\aa_d-1}' = v_i$ for $i>d$. 
The labels of all original edges are preserved. 
In the resulting graph $G[\pi]$, the set $X$ is understood to be the set of indices of 
edges incident to $v_{n+\aa_d-1}'=v_n$.

\smallskip

\begin{lemma}[Vertex splitting]\label{lem:N-ball-rec}
	Assume that 
	no loops are incident to $v_d$. 
		Then, for any  $S \subseteq X$:	
	\[   \aM_{G,\ba}(S) \ = \ \frac{1}{\aa_d!} \. \sum_{\pi \in \PM(Y_d)}  \aM_{G[\pi], \ba'}(S), \quad \text{ where } \quad  \ba' \ := \ (\aa_1,\ldots, \aa_{d-1},\underbrace{1,\ldots, 1}_{\aa_d-1}) \in \nn^{d+\aa_d-1}. \]
\end{lemma}

\smallskip

\begin{proof}
We argue by double counting the set of pairs $(O,\pi)$, 
where $O$ is an $\ba$-admissible orientation of $G$, and 
$\pi \in \PM(Y_d)$ is an involution that pairs each incoming 
edge at $v_d$ in $O$ with an outgoing edge at $v_d$, and vice versa.  
Clearly, the number of such pairs is  equal to
\[
\aa_d! \, \aM_{G,\ba}(S).
\]
Now fix such a pair $(O,\pi)$. 
We construct an orientation $O'$ of $G[\pi]$ as follows:
\begin{itemize}
	\item If $e$ is an \emph{old edge}, i.e., an edge of the original graph $G$ not incident to $v_d$, 
	then in $O'$ we orient $e$ exactly as it is oriented in $O$.
	\item If $e'$ is a \emph{new edge} of $G[\pi]$ corresponding to an edge $e$ of $G$ incident to $v_d$, 
	then in $O'$ we orient $e'$ consistently with the orientation of $e$ in $O$.
\end{itemize}
It is straightforward to check that $O'$ is $\ba'$-admissible in $G[\pi]$, 
and moreover that the pair $(O,\pi)$ can be uniquely recovered from $(\pi,O')$.  
Thus, by double counting, we obtain
\[
\aa_d! \, \aM_{G,\ba}(S) \;=\; 
\sum_{\pi \in \PM(Y_d)} \aM_{G[\pi],\ba'}(S),
\]
as claimed.
\end{proof}
\smallskip

For the remaining recurrence relations, we additionally assume that $\aa_d = 1$.  
In particular, this implies $\deg_G(v_d) = 2$, and let $e,f$ denote the two edges incident to $v_d$.  
By relabeling the edges if necessary, we will asumme that $e=e_{m-1}$ and $f=e_m$ are last two edges of $G$.
We classify $v_d$ according to the other endpoints of $e$ and $f$: it is of \defnb{type~1} if the endpoints are distinct; of \defnb{type~2} if they coincide and the common vertex lies in $\{v_1,\ldots,v_d\}$; of \defnb{type~3} if they coincide and the common vertex lies in $\{v_{d+1},\ldots,v_{n-1}\}$; and of \defnb{type~4} if they coincide and the common vertex is $v_n$.
Recall the labeling rules for  $G-e$ and  $G/e$ as described in \S\ref{subsec:rec-loop}. 

\smallskip


\smallskip

\begin{lemma}[Type 1 recurrences]\label{lem:N-ball-type-1}
Suppose that $v_d$ is of type~1.  Then, for any $S \subseteq X$:
\[\aM_{G,\ba}(S) \ = \  \aM_{G/f,\ba'}(S), \quad \text{ where } \quad \ba' \ := \  (\aa_1,\ldots, \aa_{d-1}) \in \nn^{d-1}.  \]
\end{lemma}

\smallskip

\begin{proof}
It is clear that each orientation $O$ counted by $\aM_{G,\ba}(S)$ 
corresponds to a unique orientation $O'$  counted by $\aM_{G,\ba}(S)$. 
Indeed, for every edge other than $f$, the orientation agrees in both $O$ 
and $O'$. For the edge $f$, we have: if $e$ is outgoing at $v_d$, then $f$ 
is incoming at $v_d$ in $O$; and if $e$ is incoming at $v_d$, then $f$ is 
outgoing at $v_d$. This establishes the lemma.
%
%
%
\end{proof}

\smallskip

\begin{lemma}[Type 2 and type 3 recurrences]\label{lem:N-ball-type-2-3}
	Suppose that $e$ and $f$ are incident to $v_{d-1}$.
	 Then, for any 
	$S \subseteq X$, 
	\[
	\aM_{G,\ba}(S) \ = \ 2 \.\aM_{G-e-f,\ba'}(S), 
	\qquad \text{where } \ 
	\ba' := (\aa_1,\ldots,\aa_{d-1}-1) \in \nn^{d-1}.
	\]
	
	Suppose instead that $e$ and $f$ are incident to $v_{d+1}$, with $d < n-1$. 
	Then, for any $S \subseteq X$, 
	\[
	\aM_{G,\ba}(S) \ = \ 2\. \aM_{G-e-f,\ba''}(S), 
	\qquad \text{where } \ 
	\ba'' := (\aa_1,\ldots,\aa_{d-1}) \in \nn^{d-1}.
	\]
\end{lemma}

\smallskip

\begin{proof}
We prove only the type~2 case, since the argument for type~3 is 
entirely analogous. We give a bijective proof of the lemma. 

Let $O$ be an orientation counted by $\aM_{G,\ba}(S)$. By 
$\ba$-admissibility, we must have either  
\[
(e : v_{d-1} \!\to v_d,\; f : v_d \!\to v_{d-1}) 
\quad \text{or} \quad
(e : v_d \!\to v_{d-1},\; f : v_{d-1} \!\to v_d).
\]
From $O$, we obtain an orientation $O'$ counted by 
$\aM_{G-e-f,\ba'}(S)$ by simply deleting the edges $e$ and $f$. 
One readily checks that this construction defines a bijection, with 
the factor of $2$ arising from the two possible orientations of 
$e$ and $f$.
\end{proof}

\smallskip

\begin{lemma}[Type~4 recurrences]\label{lem:N-ball-type-4}
	Suppose that $e$ and $f$ are incident to $v_{n}$, with 
	$e = e_{m-1}$ and $f = e_m$. Then, for any $S \subseteq X$, 
	\[
	\aM_{G,\ba}(S) \ = \ 
	\begin{cases}
		\aM_{G-e-f,\ba'}(S \setminus \{m-1,m\}), & \text{if } |S \cap \{m-1,m\}| = 1, \\[6pt]
		0, & \text{otherwise.}
	\end{cases},
	\]
	where \. $\ba':=(\aa_1,\ldots, \aa_{d-1}) \in \nn^{d-1}$\..
\end{lemma}

\smallskip

\begin{proof}
If $|S \cap \{m-1,m\}| \in \{0,2\}$, then $v_d$ has either no outgoing 
edges or no incoming edges, so $\aM_{G,\aa}(S) = 0$, and the claim is 
immediate. Thus we may assume $|S \cap \{m-1,m\}| = 1$.
Let $O$ be an orientation counted by $\aM_{G,\ba}(S)$. From $O$, we 
obtain an orientation $O'$ counted by $\aM_{G-e-f,\ba'}(S)$ by simply 
deleting the edges $e$ and $f$. This defines a bijection, which completes 
the proof.
\end{proof}

\medskip

\subsection{Proof of Lemma~\ref{lem:orient-ball}}
We prove the lemma by induction on the quantity 
$\aa_1 + \cdots + \aa_d$. 

If $\aa_1 + \cdots + \aa_d = 0$, then necessarily $d=0$. In this case, 
$\aM_{G,\ba}(S)$ is a constant function independent of $S$. Moreover, 
$H_X$ is a regular bipartite graph in which every vertex has degree 
$\tfrac{|X|+1}{2}$. By Hall's marriage theorem, $H_X$ admits a perfect 
matching, which implies~\eqref{eq:NMP-b} for $(H_X,\aM_{G,\ba})$.

Now suppose $\aa_1 + \cdots + \aa_d \geq 1$, which in particular forces 
$d \geq 1$. By Lemma~\ref{lem:N-ball-loop}, we may assume without loss of 
generality that $v_d$ has no loops. Applying Lemma~\ref{lem:N-ball-rec}, 
and using the fact that the sum of nonnegative functions satisfying 
\eqref{eq:NMP-b} also satisfies \eqref{eq:NMP-b}, reduces the problem to 
the case $\aa_d = 1$. (Here the value of $d$ may increase, but this is 
harmless since $\aa_1 + \cdots + \aa_d$ remains unchanged.)
There are four possibilities, according to the type of $v_d$. Specifically, 
we apply Lemma~\ref{lem:N-ball-type-1} if $v_d$ is of type~1, 
Lemma~\ref{lem:N-ball-type-2-3} if it is of type~2 or~3, 
and Lemma~\ref{lem:N-ball-type-4} if it is of type~4. In each step the 
value of $\aa_1 + \cdots + \aa_d$ strictly decreases, so the proof follows 
by induction. \qed

%
%
%
%
%
%
%
%
%
%
%
%
%

\medskip

\section{Proof of Propositions~\ref{prop:H-occ},~\ref{prop:H-ball} and Theorems~\ref{thm:NMP-occ},~\ref{thm:NMP-ball}}\label{s:proof-H-NMP}

\subsection{Proof of Proposition~\ref{prop:H-occ}}\label{ss:proof-H-occ}
Let $\sigma, \sigma'$ be arbitrary functions from $[m] \to [n]$ 
(equivalently, assignments of $m$ balls into $n$ urns). We define the 
corresponding graph $G_{\sigma,\sigma'}$ and orientation 
$O_{\sigma,\sigma'}$ as follows. The graph $G_{\sigma,\sigma'}$ has 
vertices $v_1,\ldots,v_n$ and edges $e_1,\ldots,e_m$, where each edge 
$e_i$ has endpoints $\sigma(i)$ and $\sigma'(i)$, and in 
$O_{\sigma,\sigma'}$ the edge $e_i$ is oriented from $\sigma(i)$ to 
$\sigma'(i)$. Note that $\sigma$ and $\sigma'$ can be recovered from 
$(G_{\sigma,\sigma'},O_{\sigma,\sigma'})$.

Denote by $p_{ij}$ the probability that the $i$th ball enters the $j$th 
urn. Under the correspondence above, for $S \subseteq X$, we have
\begin{align*}
	\nu_d^{\occ}(S)\,\nu_d^{\occ}(X\setminus S) 
\ 	&\propto \  \sum_{\sigma,\sigma'} \prod_{i=1}^m p_{i,\sigma(i)}\,p_{i,\sigma'(i)} \;
	\mathbf{1}\!\left\{
	\begin{array}{l}
		|\sigma^{-1}(j)|,\,|\sigma'^{-1}(j)| \geq 1 \text{ for all } j \in [d], \\[4pt]
		\sigma^{-1}(n)=S,\;\; \sigma'^{-1}(n)=X\setminus S
	\end{array}\right\},
\end{align*}
where the sum is over all functions $\sigma,\sigma':[m]\to[n]$.
This expression can be rewritten as
\[
\nu_d^{\occ}(S)\,\nu_d^{\occ}(X\setminus S) 
= \sum_{G} p_G \, \aM_{G,d}^{\occ}(S), 
\]
where the sum ranges over all graphs $G$ with $n$ vertices and $m$ edges, 
and
\[
p_G := \prod_{i=1}^m p_{i,c_i}\,p_{i,c_i'},
\]
with $v_{c_i}$ and $v_{c_i'}$ the endpoints of $e_i$ in $G$.
Proposition~\ref{prop:H-occ} now follows from the fact that 
$\aM_{G,d}^{\occ}$ satisfies~\eqref{eq:NMP-b} 
(Lemma~\ref{lem:orient-occ}), together with the observation that 
\eqref{eq:NMP-b} is preserved under positive summation.
\qed

\medskip

\subsection{Proof of Proposition~\ref{prop:H-ball}}\label{ss:proof-H-ball}
This proof follows the same argument as the proof of 
Proposition~\ref{prop:H-occ}, with $\nu_d^{\occ}$ replaced by $\nu_{\ba}$, 
$\aM_{G,d}^{\occ}$ replaced by $\aM_{G,\ba}$, and 
Lemma~\ref{lem:orient-occ} replaced by Lemma~\ref{lem:orient-ball}. 
We omit the details for brevity.
	\qed

\medskip

\subsection{Proof of Theorem~\ref{thm:NMP-occ}}\label{ss:proof-NMP-occ}
It suffices to show that for every increasing event 
$\cA \subseteq 2^{[m]}$ and every $k \geq 0$, 
\begin{align*}
	\nu_d^{\occ}(\cA, B_n=k)\,\nu_d^{\occ}(B_n=k+1) 
	\ \leq\ 
	\nu_d^{\occ}(\cA, B_n=k+1)\,\nu_d^{\occ}(B_n=k).
\end{align*}
This inequality is equivalent to 
\begin{align*}
	\sum_{S \in \binom{[m]}{k}} \ \sum_{T \in \binom{[m]}{k+1}} 
	\mathbf{1}\{S \in \cA\}\,\nu_d^{\occ}(S)\,\nu_d^{\occ}(T) 
	\ \leq\ 
	\sum_{S \in \binom{[m]}{k}} \ \sum_{T \in \binom{[m]}{k+1}} 
	\mathbf{1}\{T \in \cA\}\,\nu_d^{\occ}(S)\,\nu_d^{\occ}(T).
\end{align*}

Writing $X = S \cup T$ and $Y = S \cap T$, the inequality above is 
equivalent to 
\[
\sum_{\substack{Y \subseteq [d], \\ 0 \leq |Y| \leq k}}
\ \sum_{\substack{X \supset Y, \\ |X| = 2k+1-|Y|}}
\ \sum_{\substack{S,T \subseteq X, \\ S \cup T = X, \\ S \cap T = Y}}
\Big( \mathbf{1}\{T \in \cA\} - \mathbf{1}\{S \in \cA\} \Big) \,
\nu_d^{\occ}(S)\,\nu_d^{\occ}(T) \ \geq \ 0.
\]

Thus it suffices to prove this inequality for specific choices of $X$ and 
$Y$. By converting to an equivalent urn model if necessary, we may 
assume without loss of generality that $Y = \varnothing$ and $|X| = 2k+1$. 
In that case it suffices to show
\[
\sum_{S \in \binom{X}{k}} 
\big( \mathbf{1}\{X \setminus S \in \cA\} - \mathbf{1}\{S \in \cA\} \big)
\,\nu_d^{\occ}(S)\,\nu_d^{\occ}(X \setminus S) \ \geq \ 0.
\]

This is equivalent to proving that
\[
g_{\nu_d^{\occ}}\!\left( \cA \cap \tbinom{X}{k+1} \right) 
\ - \ 
g_{\nu_d^{\occ}}\!\left( \cA \cap \tbinom{X}{k} \right) 
\ \geq \ 0.
\]

Since $(H_{X}, g_{\nu_d^{\occ}})$ satisfies~\eqref{eq:NMP-b} by 
Proposition~\ref{prop:H-occ}, the desired inequality follows from 
\eqref{eq:HMC} (which is equivalent to~\eqref{eq:NMP-b}) together with 
the assumption that $\cA$ is an increasing event.
\qed

\medskip

\subsection{Proof of Theorem~\ref{thm:NMP-ball}}\label{ss:proof-NMP-ball}
This proof follows the same argument as Theorem~\ref{thm:NMP-occ}, with 
$\nu_d^{\occ}$ replaced by $\nu_{\ba}$ and 
Proposition~\ref{prop:H-occ} replaced by Proposition~\ref{prop:H-ball}. 
We omit the details for brevity. \qed

\bigskip

\section{Proof of Theorem~\ref{thm:CNA-occ} and Theorem~\ref{thm:CNA-ball}}\label{sec:FM-argument}

\subsection{Proof of Theorem~\ref{thm:CNA-occ}}\label{subsec:main-thm-occ}
We will use  two known results in the proof of Theorem~\ref{thm:CNA-occ}. 
The first is a theorem of Kahn and Neiman on generalized independent urn measures.

\smallskip

\begin{thm}[{\cite[Thm~14]{KN}}]\label{thm:CNC}
	The generalized  urn measures $\mu$ and $\mu^{\occ}$  satisfy \eqref{eq:CNC}.
\end{thm}

\smallskip

The second is a classical  result of Feder and Mihail~\cite{FM}.
Let $\mu$ be a probability measure on $\{0,1\}^n$. We say that $\mu$ 
satisfies the \defnb{Feder--Mihail property} if
\begin{equation*}\label{eq:FM}\tag{FM}
	\text{for every increasing $\cA \subseteq \{0,1\}^n$, there exists 
		$j \in [n]$ such that $\{x_j=1\} \uparrow \cA$,}
\end{equation*}
where $(x_1,\ldots,x_n) \in \{0,1\}^n$ is sampled according to $\mu$. 
We say that $\mu$ satisfies the \defnb{conditional Feder--Mihail property} 
if
\begin{equation*}\label{eq:CFM}\tag{CFM}
	\begin{array}{l}
			\text{every conditional measure of $\mu$, obtained by}\\ 
			\text{fixing some variables, also satisfies~\eqref{eq:FM}}.
	\end{array}
\end{equation*}

\smallskip

\begin{thm}[{\cite[Lem~3.2]{FM}}]\label{thm:Pem}
	If a measure $\mu$ on $\{0,1\}^n$ satisfies~\eqref{eq:CNC} and 
	\eqref{eq:CFM}, then $\mu$ also satisfies~\eqref{eq:CNA}.
\end{thm}

\smallskip

For references to proofs of Theorem~\ref{thm:Pem}, see \S\ref{ss:FM}.

\smallskip

%
%
%
%
%

Now let $\mu^{\occ}$ be the measure from Theorem~\ref{thm:CNA-occ}. 
By relabeling and discarding urns if necessary, we may assume without loss 
of generality that the conditioning event in $\mu^{\occ}$ is 
\[
\{\, B_1 = \cdots = B_d = 1 \,\},
\]
namely, the event that the first $d$ urns are occupied, for some $d < n$.

We now construct a new urn model from the given generalized independent urn model by 
replacing the old $n$ urns $U_1,\ldots,U_n$ with 
\. $n' := d + m(n-d)$ \.
new urns $U_1',\ldots,U_{n'}'$, according to the following deterministic 
rule:
\begin{itemize}
	\item For $j \leq d$, any balls placed in the old $j$-th urn are moved to 
	the new $j$-th urn. 
	\item For $j \geq d+1$, if the $i$-th ball is placed in the old $j$-th urn, 
	then it is moved to the new $(d + m(j-d-1) + i)$-th urn. 
\end{itemize}
In words, the new process refines each old urn $U_j$ with $j>d$ into $m$ 
sub-urns, one for each ball, so that at most one ball can land in each of 
these sub-urns.
 In particular, this implies that \. $B_j' = {B'_j}^{\occ} \in \{0,1\}$ \.  for all $j > d$.
%
Moreover, the original occupancies are recovered from 
the new ones by
\begin{align}\label{eq:old-to-new}
	B_j \ = \  
	\begin{cases}
		B_j', & j \leq d,\\
		B'_{d+m(j-d-1)+1} + \cdots + B'_{d+m(j-d-1)+m}, & j > d.
	\end{cases}
\end{align}

Let ${\mu'}^{\occ}$ denotes the law of 
$({B'_1}^{\occ},\ldots,{B'_{n'}}^{\occ})$, i.e., the occupation measure 
for the new generalized independent urn model. By the correspondence in 
\eqref{eq:old-to-new}, every increasing event 
$\cA \subseteq \{0,1\}^n$ under $\mu^{\occ}$ corresponds to a unique 
increasing event $\cA' \subseteq \{0,1\}^{n'}$ under ${\mu'}^{\occ}$. 
Hence, to prove Theorem~\ref{thm:CNA-occ}, it suffices to show that 
${\mu'}^{\occ}$, and therefore also $\mu^{\occ}$, satisfies~\eqref{eq:NA}.

Since the new model is again a generalized independent urn model, 
Theorem~\ref{thm:CNC} implies that ${\mu'}^{\occ}$ satisfies 
\eqref{eq:CNC}. By Theorem~\ref{thm:Pem}, it therefore remains only to 
verify that ${\mu'}^{\occ}$ satisfies~\eqref{eq:CFM}, which will complete 
the proof.
This is done in the following lemma, whose proof follows a standard argument in the literature showing how \eqref{eq:NMP} implies \eqref{eq:FM}; see e.g. \cite[\S3.3]{Pem} and \cite[\S2]{KN2}.


\smallskip

\begin{lemma}\label{lem:CFM}
	The measure \. ${\mu'}^{\occ}$ \.  satisfies \eqref{eq:CFM}. 
\end{lemma}

\smallskip

\begin{proof}
	First observe that for $j>d$, the event $\{B_j'=1\}$ (resp. $\{B_j'=0\}$) 
	is equivalent to the event that the $i$-th ball does (resp. does not) enter 
	the $j$-th urn, where $i \in [m]$ is determined by $i \equiv j-d \pmod m$. 
	Thus, by replacing the model with an equivalent urn model if necessary, 
	we may assume without loss of generality that there is no additional 
	conditioning on the measure, and it suffices to show that ${\mu'}^{\occ}$ 
	satisfies~\eqref{eq:FM}.
	
	Let 
	\[
	Z \ := \ B_{d+1}' + \cdots + B_{n'}'
	\]
	denotes the number of balls that entered the last $n'-d$ urns, and let 
	$\cA \subseteq \{0,1\}^{n'}$ be an arbitrary increasing event.
	 By the 
	normalized matching property~\eqref{eq:NMP} in Theorem~\ref{thm:NMP-occ} 
	(applied to the conditional urn measure where the last $n'-d$ urns are 
	merged into a single urn), we have that $Z$ is positively correlated with $\cA$. By 
	linearity of expectation, it follows that there exists some 
	$j \in \{d+1,\ldots,n'\}$ such that the event $\{B_j'=1\}$ is positively 
	correlated with $\cA$. This establishes~\eqref{eq:FM} for 
	${\mu'}^{\occ}$, completing the proof.	
\end{proof}

\smallskip

This completes the proof of Theorem~\ref{thm:CNA-occ}.\qed

\smallskip


\smallskip

\subsection{Proof of Theorem~\ref{thm:CNA-ball}}\label{subsec:main-thm-ball}
This proof follows the same argument as that of Theorem~\ref{thm:CNA-occ}, 
with the measure $\mu^{\occ}$ replaced by $\mu$ conditioned on the event 
$\{B_1 = a_1, \ldots, B_d = a_d\}$ for arbitrary $0 \leq d < n$ and 
$(\aa_1,\ldots,\aa_d) \in \nn^d$, and with 
Theorem~\ref{thm:NMP-occ} replaced by Theorem~\ref{thm:NMP-ball}. 
We omit the details for brevity. \qed

\medskip

\medskip

\section{Final remarks} \label{s:finrem}

\subsection{Generalized  interval urn measures}
The {generalized independent urn occupation measures} and 
{enumeration measures} are special cases of a broader 
class called \defnb{generalized independent interval urn measures}.
For each $j \in [n]$, let us be given a sequence of cutpoints
\. 
$0 = c_0(j) < c_1(j) < \cdots < c_{k_j}(j) = m+1$\..
For an urn assignment $\sigma:[m] \to [n]$, define
\[
\ax_j(\sigma) \;=\; t 
\quad \text{iff} \quad 
c_t(j) \,\leq\, B_j \,<\, c_{t+1}(j),
\]
that is, $\ax_j(\sigma)$ is  the index of the interval 
containing  the number of 
balls in the $j$-th urn.
We write $\mu^{\mathrm{int}}$ for the probability measure of the 
random vector $(\ax_1,\ldots,\ax_n) \in \mathbb{N}^n$.

The measure is called a \defnb{generalized threshold urn measure} when 
$c_{2}(j) = \infty$ for all $j \in [n]$, since each urn is classified 
according to whether its ball count exceeds the threshold $c_1(j)$. 
In this case the measure is supported on $\{0,1\}^n$.

\smallskip

The following question is a natural generalization of Theorems~\ref{thm:CNA-occ}  and \ref{thm:CNA-ball}.

\smallskip

\begin{question}[{\cite[\S5]{KN}}]\label{quest:interval}
	Are generalized threshold (resp. interval) measures \eqref{eq:CNA}?
\end{question}

\smallskip
We remark that the argument of Theorem~\ref{thm:CNA-ball} can be adapted 
to the special case of interval measures in which 
\. $c_{t+1}(j) - c_t(j) \in \{1,2\}$ \. for all \. $0 \leq t \leq m$ \. and 
\. $1 \leq j \leq n$. The proof follows by the same reasoning, and we omit it for conciseness.
Beyond this special 
case, if the answer to Question~\ref{quest:interval} is affirmative, its 
proof may well require new ideas beyond those developed in this paper.
On the other hand, a negative answer to Question~\ref{quest:interval} in 
the case of threshold measures would be especially interesting, as it 
would disprove Pemantle’s conjecture that \eqref{eq:CNC} implies 
\eqref{eq:CNA}.

\smallskip

%

\subsection{Feder--Mihail arguments}\label{ss:FM}
The classical Feder--Mihail Theorem~\ref{thm:Pem} traces back to \cite[Lem.~3.2]{FM}, where it was 
stated for balanced matroids. The proof was presented in the special case 
where one of the monotone events is a single-coordinate event, i.e., of 
the form $\{\ax_i = 1\}$, and the authors remarked that the argument 
extends to arbitrary monotone events. For a proof under the assumptions 
used in this paper, see \cite[Thm.~1.3]{Pem}, although that proof likewise 
presents only the case where one monotone event is of the form 
$\{\ax_i = 1\}$. A proof for arbitrary monotone events can be found in 
\cite[Ex.~4.6]{LP}.

\smallskip

\subsection{Log-concavity}\label{ss:lc}
For a measure 
$\mu$ on $\{0,1\}^n$, the \defnb{rank sequence} is 
$(r_j)_{j=0}^n$, where
\[
r_j \ := \ \mu\big\{ (\ax_1,\ldots,\ax_n) \in \{0,1\}^n \ \mid \ 
\ax_1+\cdots+\ax_n = j \big\}.
\]
We say that $\mu$ is \defnb{ultra-log-concave} if its rank sequence has no 
internal zeros and satisfies
\begin{equation*}\tag{ULC}\label{eq:ULC}
	\left(\frac{r_j}{\binom{n}{j}}\right)^2 
	\ \geq \ 
	\frac{r_{j+1}}{\binom{n}{j+1}}
	\cdot 
	\frac{r_{j-1}}{\binom{n}{j-1}}
	\quad \text{for } 0<j<n.
\end{equation*}
The interplay between negative correlation and log-concavity has been 
investigated for several decades and remains an active and expanding area 
of research; see \cite{Sta,Bre,Pem,BBL,KN2,BH,HSW,CKN} and the works cited therein 
for a non-exhaustive list.

Of particular relevance to conditional measures is a result of Pemantle, 
which shows that \eqref{eq:CNC}, \eqref{eq:CNA}, and \eqref{eq:ULC} are 
equivalent for exchangeable measures~\cite[Thm.~2.7]{Pem}. 
In general, however, there are no logical implications between  \eqref{eq:ULC} and \eqref{eq:CNC} (resp. \eqref{eq:CNA}).

The following gives a natural example showing that \eqref{eq:ULC} does not 
imply \eqref{eq:CNC} (and hence \eqref{eq:CNA}). Consider the uniform measure on 
independent sets of a matroid. Verifying that this measure satisfies 
\eqref{eq:ULC} was a long-standing problem of Mason~\cite{Mas}, resolved independently 
by Anari et al.~\cite{ALOV} and Br\"and\'en--Huh~\cite{BH}.
In contrast, a counterexample to \eqref{eq:CNC} arises from any matroid in 
which there exist two elements that are positively correlated under the uniform 
distribution on the bases of the matroid. The earliest construction of such a 
matroid appears in the work of Seymour and Welsh~\cite[p.~495]{SW}; see 
also~\cite{HSW} for further examples. Now, modify such a matroid by adding $k$ 
new parallel elements to each element of its ground set. As $k \to \infty$, the 
uniform measure on independent sets concentrates on the bases, which implies 
that this measure fails to satisfy \eqref{eq:NC} (and hence \eqref{eq:CNC}) for 
sufficiently large $k$.\footnote{This construction is due  to Petter Br\"and\'en, personal communication, September 2nd 2025.}

Conversely, earliest examples showing that \eqref{eq:CNA} (and hence \eqref{eq:CNC}) 
does not imply \eqref{eq:ULC} can be found in~\cite[\S7]{BBL} and 
\cite[Thm.~6]{KN2}. Interestingly, the urn occupation measures $\mu^{\occ}$ studied in this 
paper provide perhaps a more natural counterexample, as even the ordinary 
occupation measure is already known to fail \eqref{eq:ULC}; 
see~\cite[Ex.~31]{KN2}.

\smallskip




\subsection{External fields}\label{ss:Rayleigh}
For a measure $\mu$ on $\{0,1\}^n$, we can \defnb{impose an external field} $\mathbf{W}=(W_1,\ldots,W_n) \in \mathbb{R}_{\geq 0}^n$ by defining
\[
{\mathbf{W}} \circ \mu (\bx) \ \propto \ \mu(\bx)\,\prod_{i=1}^n W_i^{\,\ax_i},
\qquad \bx \in \{0,1\}^n.
\]
We say that $\mu$ satisfies ({\color{blue}{NA+}})
 if \. ${\mathbf{W}} \circ \mu $ \. satisfies 
 \eqref{eq:NA}
for every external field $\mathbf{W} \in \rr_{\geq 0}^n$.
Note that ({\color{blue}{NA$+$}})  implies \eqref{eq:CNA}, by taking the  limits \. $\mathrm{W}_i \to \infty$ \. or  \. $\mathrm{W}_i \to 0$\., corresponding to  conditioning on  \. $\{\ax_i=1\}$ \. or  \. $\{\ax_i=0\}$\., respectively.

The urn occupation measure $\mu^{{occ}}$ considered in this paper 
provides a natural example of a measure that satisfies \eqref{eq:CNA}
but not ({\color{blue}{NA$+$}}), 
already in the ordinary case; see \cite[Ex.~32]{KN2}.
The same phenomenon occurs  if we replace  \eqref{eq:CNA}, ({\color{blue}{NA+}}) with  \eqref{eq:CNC},  ({\color{blue}{NC+}}).
In the literature, ({\color{blue}{NC$+$}}) is also referred to as the \defnb{Rayleigh property}.

\smallskip

\subsection{Stochastic covering property}\label{ss:SCP}
The normalized matching property \eqref{eq:NMP} in \S\ref{ss:NMP} is closely related to another property from the negative dependence literature, known as the  \defnb{stochastic covering property}, due to Pemantle and Peres~\cite[Def.~2.1]{PP}.
Let $\nu$ be a measure on $\{0,1\}^{m}$.
We say that $\nu$ satisfies \eqref{eq:SCP} if, for all  \. $\ba, \mathbf{b} \in \{0,1\}^m$ \. satisfying $\ba \leq \mathbf{b}$   and all \. $S \subseteq [m]$: 
\begin{equation*}\tag{SCP}\label{eq:SCP}
	\nu(\cA \mid \az_i=\aa_i \text{ for } i \in S) \  \leq \ \nu(\cA \mid \az_i=\mathrm{b}_i \text{ for } i \in S) \quad \text{ for increasing } \cA \subseteq 2^{[m]},
\end{equation*}
where  $\bZ:=(\az_1,\ldots, \az_m)$ is  sampled according to $\nu$.
This is closely analogous to~\eqref{eq:NMP}, except that here the conditioning 
is on the coordinates $(\az_i)_{i \in S}$, rather than on the sum $\az_1+\ldots+\az_m$.


Interestingly, the measure $\nu_d^{\occ}$ from \S\ref{ss:NMP} 
(recall that this is the distribution of the set of balls entering the $n$-th urn, 
conditioned on the first $d$ urns being non-empty) provides a natural example 
of a measure that satisfies~\eqref{eq:NMP} but not~\eqref{eq:SCP}. 
Indeed, the measure $\nu_d^{\occ}$ satisfies~\eqref{eq:NMP}, 
as proved in Theorem~\ref{thm:NMP-occ}. 
On the other hand, the following example shows that this measure does not satisfy~\eqref{eq:SCP}. 

Consider the generalized independent urn model with $m=2$ balls and $n=2$ urns, 
where each ball has equal probability of entering either urn. 
Let $d=1$, and let $\cA$ be the event that the second ball enters the second urn. 
Then we have
\[
\nu_d^{\occ}(\cA \mid \az_1=0) \ = \ \tfrac{1}{2}; 
\qquad    
\nu_d^{\occ}(\cA \mid \az_1=1) \ = \ 0. 
\]

For the first equation, note that $\az_1=0$ means the first ball does not enter the second urn, 
so it must enter the first urn, which ensures that the first urn is occupied. 
In this case, the second ball then has equal chance to enter either urn, so 
\. $\nu_d^{\occ}(\cA \mid \az_1=0) = \tfrac{1}{2}$\.. 
For the second equation, $\az_1=1$ means the first ball enters the second urn. 
To keep the first urn occupied, the second ball must then enter the first urn, 
which implies 
\. $\nu_d^{\occ}(\cA \mid \az_1=1) = 0$. 
This demonstrates that $\nu_d^{\occ}$ does not satisfy~\eqref{eq:SCP}.

\subsection{Admissible orientations}
The approach of employing admissible orientations to study urn problems was first introduced by Kahn and Neiman in \cite[\S4]{KN}, and the proofs of Lemmas~\ref{lem:orient-occ} and \ref{lem:orient-ball} are guided by their approach. It is therefore natural to speculate that admissible orientations may also play an important role in resolving Question~\ref{quest:interval}.

\smallskip

%
%
%
%



\vskip.7cm

\subsection*{Acknowledgements}
The author is very  grateful to Jeff Kahn for introducing the  problems to the author   and for many stimulating discussions.
The author would also like to thank  Petter Br\"{a}nd\'{e}n, Igor Pak and Robin Pemantle for valuable references 
and suggestions,  and   Milan Haiman, Bhargav Narayanan,  Sam Spiro, and Natalya Ter-Saakov for inspiring discussions. This research was partially supported by NSF grant DMS--2246845.

\vskip.9cm


\end{document}